\documentclass[JEDP]{cedram-sem}
\usepackage{amsfonts}
\usepackage{amsthm,mathtools,bm,dsfont,amssymb}

\newtheorem{theorem}{Theorem}[section]

\newtheorem{corollary}[theorem]{Corollary}

\newtheorem*{remarks*}{Remarks}
\newtheorem*{remark*}{Remark}

\newcommand{\pt}{\partial}
\newcommand{\uu}{\mathbf{u}}
\newcommand{\QQ}{\mathbf{Q}}
\newcommand{\Dh}{|\nabla|}
\newcommand{\R}{\mathbb{R}}
\newcommand{\Sb}{\mathbb{S}}
\newcommand{\T}{\mathbb{T}}
\newcommand{\C}{\mathbb{C}}
\newcommand{\eps}{\varepsilon}
\newcommand{\im}{{\mathrm{i}}}
\newcommand{\Xsp}{\mathbf{X}}
\newcommand{\Usp}{\mathbf{U}}

\title 
[Short Primer on the Half-Wave Maps Equation]
{A Short Primer on the Half-Wave Maps Equation}

\author
[E. \lastname{Lenzmann}]
{\firstname{Enno} \lastname{Lenzmann}}
\address{University of Basel \\ Department of Mathematics \\ Spiegelgasse 1, CH-4051 Basel, Switzerland}
\thanks{E.~L. is supported by the Swiss National Science Foundation (SNSF) under grant no.~200021-149233.}
\email{enno.lenzmann@unibas.ch}

\keywords{}

\begin{document}

\begin{abstract}
     We review the current state of results about the half-wave maps equation on the domain $\R^d$ with target $\Sb^2$. In particular, we focus on the energy-critical case $d=1$, where we discuss the classification of traveling solitary waves and  a Lax pair structure together with its implications (e.\,g.~invariance of rational solutions and infinitely many conservation laws on a scale of homogeneous Besov spaces). Furthermore, we also comment on the one-dimensional space-periodic case. Finally, we list some open problem for future research.
     \end{abstract}

%

\maketitle

\section{Introduction}

This expository note is intended to give an overview on results about the \textbf{half-wave maps equation} posed in $\R^d$ with target $\Sb^2$ (embedded in $\R^3$). The corresponding equation reads
\begin{equation} \tag{HWM}
\pt_t \uu = \uu \times \Dh \uu 
\end{equation}
for the function $\uu : [0,T) \times \R^d \to \Sb^2$. Here $\times$ denotes the standard vector product in $\R^3$ and the operator $|\nabla|$ is defined via the Fourier transform as $\widehat{(|\nabla| \uu})(\xi) =  |\xi| \widehat{\uu}(\xi)$. Below we also address the one-dimensional periodic setting when $\R$ is replaced by  $\T$.

In what follows, we will mainly deal with $d=1$ space dimensions. As a matter of fact, this case displays a very rich list of interesting analytical phenomena such as {\em energy-criticality,  conformal invariance,} and a {\em Lax pair structure}. We try to highlight the relevant mathematical results that have been found so far, followed by an outline of future open problems. 

Let us also mention that, from a physical perspective, the one-dimensional half-wave maps equation is also of significant interest, since it (formally) arises as a semi-classical and continuum limit of {\em Haldane-Shastry (HS) spin chains} and {\em classical spin systems of Calogero-Moser (CM) type}. Models of (HS) and (CM) are exactly solvable and completely integrable quantum and classical systems, respectively. For more details on this, we refer the reader to the appendix of \cite{LeSc-18}.

\section{Some Basics Facts}

We start with collecting some fundamental properties of the geometric evolution equation given by (HWM).

\subsection{Conservation Laws and Hamiltonian Structure}

A moment's reflection shows that (HWM) exhibits (formal) conservation of the energy 
\begin{equation}
E[\uu] = \frac 1 2 \int_{\R^d} \uu \cdot |\nabla| \uu \, dx = c_d \iint_{\R^d \times \R^d} \frac{|\uu(x)-\uu(y)|^2}{|x-y|^2} \, dx \, dy,
\end{equation}
with some constant $c_d > 0$ (where the last identity is a classical fact from harmonic analysis). From this observation, we see that the corresponding energy space is found to be the homogeneous Sobolev space
$$
\dot{H}^{\frac 1 2}(\R^d; \Sb^2) = \{ \uu : \R^d \to \R^3 \mid \mbox{$E[\uu] < +\infty$ and $|\uu| = 1$ a.\,e.} \}.
$$
It is not hard to see that (HWM) can indeed be written as a Hamiltonian equation of motion $\pt_t \uu = \{ \uu, E \}$ with the canonical Poisson bracket for $\Sb^2$-valued functions defined as
\begin{equation}
\{ u_i(x), u_j(y) \} = \eps_{ijk} u_k(x) \delta(x-y),
\end{equation}
where $\eps_{ijk}$ is the anti-symmetric Levi-Civit\`a symbol.  The energy space $\dot{H}^{\frac 1 2}(\R^d; \Sb^2)$ can then be seen as a phase space for the infinite-dimensional Hamiltonian system whose equation of motions are given by the half-wave maps equation.

Another type of (formally) conserved quantities for (HWM) are found to be
\begin{equation}
\mathbf{S}[\uu] = \int_{\R^d} (\uu - \mathbf{P}) \, dx, \quad M[\uu] = \int_{\R^d} |\uu - \mathbf{P}|^2 \, dx,
\end{equation}
where we require that $\uu - \mathbf{P} \in L^1(\R^d)$ or $\uu - \mathbf{P} \in L^2(\R^d)$, respectively, for some constant  point $\mathbf{P}$ on $\Sb^2$. In physical terms, the quantity $\mathbf{S}[\uu]$ is the total spin of the system represented by the field $\uu$. Of course, the rotational symmetry on the domain $\R^d$ (i.\,e.~the continuous action of the group $\mathrm{SO}(d)$ on the domain $\R^d$ for $d \geq 2$) also induces formal conservation laws by Noether's theorem corresponding to the conservation of {\em angular momentum}. However, we omit any further details here because the conservation of angular momentum does not play any role in what follows. Another -- and more interesting conservation law with respect to analysis -- is related to the invariance under spatial translations in the domain $\R^d$. For example when $d=1$, we obtain the (formal) conversation of the quantity 
\begin{equation}
P[\uu] = \int_{\R} \frac{u_2 \pt_x u_1 - u_1 \pt_x u_2}{1-u_3} \, dx , 
\end{equation}
which can be regarded as a {\em linear momentum}. For more details and the geometric meaning of $P[\uu]$, we refer to \cite{LeSc-18}.

With regard to conserved quantities, we finally mention that the one-dimensional (HWM) in fact possesses {\em infinitely many conservation laws} due to the existence of a {\em Lax pair}; see Section \ref{sec:Lax} for more on this remarkable feature that indicate some kind of complete integrability in the one-dimensional setting.

\subsection{Criticality and Conformal Symmetry}

For a constant $\lambda >0$, we readily check that rescaling 
$$
\uu(t,x) \mapsto \uu_\lambda(t,x) = \uu(\lambda t, \lambda x)
$$ 
maps solutions of (HWM) into solutions defined on the time interval $[0, \lambda^{-1} T)$. Likewise, the energy transforms as 
$$
E[\uu_\lambda] = \lambda^{d-1} E[\uu].
$$
Hence, we see that $E[\uu_\lambda] = E[\uu]$ in $d=1$ space dimensions, which implies that the evolution problem is {\em energy-critical} in one space dimension. In higher dimensions $d \geq 2$, the equation (HWM) becomes {\em energy-supercritical}.

Another feature of the energy-critical case when $d=1$ is the {\em conformal invariance} of the energy. More precisely, let $\uu \in \dot{H}^{\frac 1 2}(\R; \Sb^2)$ be given and consider its harmonic extension $\uu^e : \R_+^2 \to \R^3$ to the upper half-plane $\R^2_+$. Then it is a classical fact that
\begin{equation}
E[\uu] = \frac 1 2 \int_{\R} \uu(x) \cdot |\nabla| \uu(x) \, dx = \frac 1 2 \iint_{\R^2_+} |\nabla_{(x,y)} \uu(x,y)|^2 \, dx \, dy. 
\end{equation} 
Now the right-hand side is invariant under conformal transformations $\phi : \R_+^2 \to \R_+^2$ preserving the upper half-plane $\R_2^+ \simeq \C_+$, which is known to be the M\"obius subgroup $\mathrm{PSL}(2, \R)$. An elementary calculation yields that
$$
E[\uu \circ \phi] = E[\uu] \quad \mbox{for all $\phi \in \mathrm{PSL}(2, \R)$}.
$$
This conformal invariance of $E[\uu]$ will play a role in the classification of traveling solitary waves (in addition to another action of the M\"obius group on the target $\Sb^2$).

\subsection{Geometric Formulation and Variants}

In analogy to the Schr\"odinger maps equation $\partial_t \uu = \uu \times \Delta \uu$ for $\Sb^2$-valued maps, we can recast the half-wave maps equation into a more geometric form by exploiting that $\Sb^2$ is a K\"ahler manifold. Indeed, we can write (HWM) as
\begin{equation} \label{eq:HWM_kaehler}
\pt_t \uu = J_\uu P_\uu \Dh \uu 
\end{equation} 
where $J_\uu = \uu \times : T_\uu \Sb^2 \to T_\uu \Sb^2$ is the standard complex structure on $\Sb^2$ and $P_\uu : \R^3 \to T_\uu \Sb^2$ denotes the projection onto the tangent space $T_\uu \Sb^2$. So far, this reformulation of (HWM) has not been proven to be fruitful in the analysis yet. However, from this point of view, a natural modification of (HWM) arises when the compact K\"ahler manifold $\Sb^2$ is replaced by the non-compact K\"ahler $\mathbb{H}^2$ (hyperbolic plane). For more details on this, see Section \ref{sec:h2} below. 

Furthermore, in analogy to the Landau-Lifshitz-Gilbert equation (which is a combination of Schr\"odinger and harmonic maps heat flow), we can generalise the equation \eqref{eq:HWM_kaehler} by adding a parabolic term, i.\,e.,
\begin{equation}
\pt_t \uu = \alpha J_\uu P_\uu \Dh u - \beta P_\uu \Dh u
\end{equation}   
with constants $\alpha, \beta \geq 0$. For $\alpha = 0$, this is {\em half-harmonic maps heat flow} with target $\Sb^2$. See \cite{SiWeZh-18b} for infinite-time blow-up solutions for the purely parabolic evolution problem (when $\alpha=0$) with target $\Sb^1$.

\subsection{Cauchy Problem}

Up to now, the Cauchy problem for (HWM) lacks a full-fledged well-posedness theory for initial data in the energy space $\dot{H}^{\frac 1 2}(\R^d; \Sb^2)$. In high dimensions $d \geq 5$ and for small energy initial data, the recent work in \cite{KrSi-18} establishes global-in-time existence. However, these techniques do not seem to be applicable for low dimensions $d \leq 3$; in particular, the energy-critical case $d=1$ seems out of scope by these techniques. 

Let us mention that short-time existence of solutions for initial data in $\dot{H}^{k}$ with $k > \frac{d}{2} + 1$ can be established by standard arguments. However, the global-in-time existence or finite-time blowup is a completely open question. Also, by standard approximation arguments (e.\,g.~parabolic regularisation), we can deduce existence of global weak solutions for finite energy data (where uniqueness of weak solutions is of course an open problem). See, e.\,g., \cite{PuGu-13} for such results on the domain $\mathbb{T}^d$ with $d=1,2,3$, but the arguments there can be carried over to domain $\R^d$ with $d=1,2,3$.

\section{Traveling Solitary Waves} \label{sec:solitary}

We will now turn out attention to the energy-critical (HWM) in $d=1$ space dimension. We seek special solutions that are given by {\em traveling solitary waves}. By definition, these solutions are of the form
\begin{equation}
\uu(t,x) = \QQ_v(x-vt),
\end{equation}
where $v \in \R$ is a given constant velocity. The profile function $\QQ_v : \R \to \Sb^2$ is seen to satisfy the nonlinear equation
\begin{equation} \label{eq:Qv}
\QQ_v \times |\nabla| \QQ_v + v \pt_x \QQ_v = 0.
\end{equation}
Indeed, we shall be interested in finite-energy solutions and hence we shall assume that $\QQ_v \in \dot{H}^{\frac 1 2}(\R; \Sb^2)$. Other non-trivial solutions $\QQ_v$ of \eqref{eq:Qv} with infinite energy are also given below.

In the special case $v=0$, corresponding to static solutions of (HWM), leads to the {\em half-harmonic maps equation}
\begin{equation} \label{eq:Qv_0}
\QQ \times \Dh \QQ = 0
\end{equation} 
for maps $\QQ \in \dot{H}^{\frac 1 2}(\R; \Sb^2)$. In fact, this equation was introduced in \cite{DaLioRi-11} as a model problem to study higher regularity for a nonlocal elliptic energy-critical problem with conformal symmetry (in analogy to harmonic maps from $\R^2$ into $\Sb^2$). See also \cite{MiSi-15}. 

\subsection{Complete Classification of Traveling Solitary Waves}

The following classification result obtained in \cite{LeSc-18} now extends the known classification result for half-harmonic maps (see e.\,g.~\cite{MiSi-15}). However, the presence of the term $v \pt_x \QQ_v$ in \eqref{eq:Qv} calls for some further ideas, since $\QQ_v$ with $v \neq 0$ do not correspond technically speaking to free boundary minimal disks but to a more general case. From \cite{LeSc-18} we recall the following result.

\begin{theorem} \label{thm:Qv_class}
Let $v \in \R$ and suppose $\QQ_v \in \dot{H}^{\frac 1 2}(\R; \Sb^2)$ solves \eqref{eq:Qv}. Then the following holds.
\begin{enumerate}
\item[(i)] If $|v| < 1$, then $\QQ_v$ must be of the form\footnote{Taking $-v$ here corresponds to $+v$ in \cite{LeSc-18}, where a different sign convention for (HWM) is used.}
$$
\QQ_v(x) = R \left ( \sqrt{1-v^2} \, \mathrm{Re} \, \mathcal{B}(x), \mp \sqrt{1-v^2} \, \mathrm{Im} \, \mathcal{B}(x), \mp v \right )
$$ 
with some fixed rotation $R \in \mathrm{SO}(3)$ and $\mathcal{B} : \overline{\C}_+ \to \C$ is a finite Blaschke product of degree $m \in \mathbb{N}_0$, i.\,e.,
$$
B(z) = \prod_{k=1}^m \frac{z-z_k}{z-z_k^*} 
$$
with arbitrary points $z_1, \ldots, z_m$ in the complex upper half-plane $\C_+$.
\item[(ii)] If $|v| \geq 1$, then $\QQ_v(x) \equiv \mathbf{P}$ for some constant $\mathbf{P} \in \Sb^2$.
\end{enumerate}
\end{theorem}

\begin{remarks*}
{\em 1. In (i), the trivial case $\QQ_v \equiv \mbox{const}.$ corresponds to degree $m=0$. 

2. The energy of the maps $\QQ_v \in \dot{H}^{\frac 1 2}(\R; \Sb^2)$ for $|v| < 1$ is found to be
$$
E[\QQ_v] = (1-v^2) \cdot m \pi.
$$
In particular, we see that $E[\QQ_v] \to 0$ as $|v| \to 1^-$. As an interesting consequence, we can construct traveling solitary waves with arbitrarily small energy. This is in stark contrast to other energy-critical evolution PDEs (e.\,g.,~energy-critical wave maps, Schr\"odinger maps, NLW and NLS equations) where small energy data lead to scattering to free solutions.

3. For degree $m=1$, the profiles $\QQ_v$ can regarded as {\bf ground states}, i.\,e. as the nontrivial solutions with smallest energy. For example if we take $z_1=+\im$ and $v=0$, then
$$
\QQ(x) = \left ( \frac{x^2-1}{x^2+1}, \frac{2x}{x^2+1}, 0 \right ) 
$$  
is a half-harmonic map, i.\,e., a finite-energy solutions of \eqref{eq:Qv_0} with $v=0$. More generally, we see that all half-harmonic maps are (up to rotations on the sphere) are given by rational parametrisations of the equator, where the integer $m \geq 1$ is the winding number.

4. In fact, one can view the map $\QQ_{v=0} \mapsto \QQ_{v \neq 0}$ as {\em M\"obius transform} (conformal transformation) on the target $\Sb^2$. See \cite{LeSc-18} for more details and why this transformation can be seen as a {\em Lorentz boost} implemented by the M\"obius group acting on $\Sb^2$. The presence of the factor $\sqrt{1-v^2}$ can be viewed as a kind of {\em Lorentz contraction} occurring for relativistic Lorentz boosts. 

5. For $|v| < 1$ and integer $m \geq 1$, the maps 
$$
\tilde{\QQ}_v(x) = (\sqrt{1-v^2} \cos (mx), \sqrt{1-v^2} \sin(mx), -v ) \in \dot{H}^{\frac 1 2}_{\mathrm{loc}}(\R; \Sb^2)
$$
are also (distributional) solutions of \eqref{eq:Qv}. However, these maps do not belong to $\dot{H}^{\frac 1 2}$ and thus the corresponding solutions $\uu(t,x)= \tilde{Q}_v(x-vt)$ are traveling waves for (HWM) with infinite energy. It seems an interesting open problem to show that $\tilde{Q}_v$ provide all infinite-energy traveling waves for (HWM) up to symmetries.

6. See also \cite{ZhSt-15} for a discussion of traveling solitary waves for the half-wave maps equation (HWM) in the physics literature.
}
\end{remarks*}

\begin{proof}[Sketch of the Proof of Theorem \ref{thm:Qv_class}] 
Proving cases (i) and (ii) are of very different nature. 

\medskip
{\em Sketching the Proof of (i).} We briefly summarise the arguments in \cite{LeSc-18} used to obtain (i). Let $|v| < 1$ be given. By a delicate bootstrap argument (building upon the regularity theory in \cite{DaLioRi-11} for half-harmonic maps), we first show that any $\QQ_v \in \dot{H}^{\frac 1 2}$ solving \eqref{eq:Qv}  belongs to $C^\infty \cap \dot{H}^2$.  For this regularity argument, the assumption $|v| < 1$ is crucial. 

Once the regularity of $\QQ_v$ is improved (actually, some higher H\"older continuity would suffice), we show as a next step that -- after a suitable rotation on $\Sb^2$ -- we have
\begin{equation} \label{eq:Qv_plane}
\QQ_v(x) = ( \sqrt{1-v^2} \, f(x),  \sqrt{1-v^2} \, g(x), \pm v ) 
\end{equation} 
with some smooth function $f, g: \R \to \R$ such that $f^2 + g^2 \equiv 1$. In geometric terms, equation \eqref{eq:Qv_plane} means that the image of $\QQ_v$ lies in a plane $E$ that is parallel to the equatorial plane $\{z=0\}$ and $E$ has distance $|v|$ to the plane $\{ z= 0 \}$. Furthermore, since we have $\QQ_v(-\infty) = \QQ_v(+\infty)$ by the finite energy condition, we see that $x \mapsto \QQ_v(x)$ parametrises a circle on $\Sb^2$ with radius $\sqrt{1-v^2}$.

But how to prove \eqref{eq:Qv_plane} above? To attack this problem, we consider $\QQ_v : \R \to \Sb^2$ as the boundary curve of a minimal surface $\Sigma$ inside the the unit ball $B_1(0) \subset \R^3$. Thanks to the finite energy assumption $\QQ_v \in \dot{H}^{\frac 1 2}$, we find that $\QQ_v(-\infty) = \QQ_v(+\infty)$ showing that $Q_v$ yields a closed curve on the unit sphere $\Sb^2$.

Next we let $\QQ_v^e : \R^2_+ \to \R^3$ denote the (unique) bounded harmonic extension of $\QQ_v$ to the upper half-plane $\R_+^2$.  By the strong maximum principle, it follows that $\QQ_v^e(\R_+^2) \subset B_1(0)$ holds, unless $\QQ_v$ is constant (which is a trivial case that we exclude here). Now, we apply the techniques of {\em Hopf differentials} to prove that the map $\QQ_v^e: \R_+^2 \to \R^3$ traces out a {\em minimal surface} $\Sigma \subset B_1(0)$ whose boundary $\pt \Sigma$ given by the parametrised curve $\QQ_v : \R \to \Sb^2$. To carry out this argument, we identify $\R_+^2$ with the complex upper half-plane $\C_+$ via $z = x+ \im y$. Next, we define the functions $\mathrm{\tt Hopf}_{\nu} : \C_+ \to \C$ with $\nu=1,2$ by setting
\begin{equation}
\mathrm{\tt Hopf}_1(z) = \pt_z \QQ_v^e \cdot \pt_z \QQ_v^e  \quad \mbox{and} \quad \mathrm{\tt Hopf}_2(z) = \pt_z^2 \QQ_v^e \cdot \pt_{z}^2 \QQ_v^e,
\end{equation}
where $\pt_z = \frac 1 2 (\pt_x - \im \pt_y)$ and $U \cdot V = \sum_{j=1}^3 U_j V_j$ for $U, V \in \C^3$. From the harmonicity $\Delta \QQ_v^e = 0$, we directly infer that $\mathrm{\tt Hopf}_1$ and $\mathrm{\tt Hopf}_2$ are both holomorphic on $\C_+$. Now by using the equation \eqref{eq:Qv} satisfied by $\QQ_v = \QQ_v^e |_{\pt \R_+^2}$, we deduce that $\mathrm{\tt Hopf}_1(z) \equiv 0$ vanishes identically. But this shows that $\QQ_v^e$ is a conformal map from $\R^2_+$ into $B_1(0)$, i.\,e., we have $|\pt_x \QQ_v^e | = |\pt_y \QQ_v^e |$ and $\pt_x \QQ_v \perp \pt_y \QQ_v^e$. By its harmonicity, this means that $\QQ_v^e$ traces out a minimal surface $\Sigma \subset B_1(0)$. Finally, we prove that the second Hopf differential $\mathrm{\tt Hopf}_2(z) \equiv 0$ also vanishes identically. From this fact we can deduce that $\Sigma$ is indeed a {\em flat disk}, which implies that its boundary $\partial \Sigma$ is a circle on $\Sb^2$. This completes the sketch of the proof of \eqref{eq:Qv_plane}.

The rest of the proof of case (i) now boils down to a problem in complex analysis. By the previous discussion, the harmonic extension $\QQ_v^e$ can be identified with a bounded holomorphic function $F : \C_+ \to C$ satisfying the boundary condition $|F|^2 \equiv 1$ on $\pt \C_+ \simeq \R \times \{0\}$. Thus $F$ is an inner function and hence has the canonical factorisation
\begin{equation}
F(z) = \lambda e^{\im \alpha z} B(z)S(z),
\end{equation}
where $\lambda \in \C$, $|\lambda| = 1$ and $\alpha \geq 0$. Here $B(z)$ is (a possibly infinite) Blaschke product with zeros on $\C_+$ and $S$ is the so-called singular inner part. However, by the regularity estimates for $F$ on the boundary $\pt \C_+$, we can show that $S(z) \equiv 1$ is trivial. Furthermore, by the finite energy property, we have 
$$
-\int_{\pt \C_+} \overline{F} \pt_y F |_{y=0} \,dx = \int_{\R} \QQ_v \cdot |\nabla| \QQ_v \, dx < +\infty, 
$$ 
we find that $\alpha = 0$ and that $B(z)$ must be a finite Blaschke product, i.\,e., we have
$$
B(z) = \prod_{k=1}^m \frac{z-z_k}{z-z_k^*}
$$
with some $z_1, \ldots, z_m \in \C_+$ and $m \in \mathbb{N}$. The interested reader may consult \cite{LeSc-18} for details.

\medskip
{\em Sketching the Proof of (ii).} In the case $|v| \geq 1$, no regularity method is known to improve the regularity of maps $\QQ_v \in \dot{H}^{\frac 1 2}$ solving \eqref{eq:Qv}. However, we can build up an `Pohozaev-type' argument to show that finite-energy solutions of \eqref{eq:Qv} must be constants whenever $|v| \geq 1$ holds. Luckily, no higher regularity is needed to carry out this argument.

Indeed, let $\QQ_v^e : \R^2_+ \to \R^3$ be the harmonic extension of $\QQ_v : \R \to \Sb^2$. Since $\QQ_v \in \dot{H}^{\frac 1 2}(\R; \Sb^2)$, we have that $\QQ_v^e \in\dot{H}^1(\R_+^2;\R^3)$. By testing the equation \eqref{eq:Qv} against the Hilbert transform\footnote{To be understood modulo constants. More precisely,  for $\QQ_v \in \dot{H}^{\frac 1 2}(\R; \Sb^2)$, there exists some constant $\mathbf{P} \in \Sb^2$ such that $\QQ_v - \mathbf{P} \in D^{\frac 1 2}(\R)$ = closure of $C^\infty_c(\R)$ with respect to $\| \cdot \|_{\dot{H}^{\frac 1 2}}$. Thus $H(\QQ_v)$ means actually $H(\QQ_v-\mathbf{P})$.} $H(\QQ_v) \in \dot{H}^{\frac 1 2}$, we find after some calculations the following identity:
\begin{equation} \label{eq:poho}
v \iint_{\R_+^2} (|\pt_x \QQ_v^e|^2 + |\pt_y \QQ_v^e|^2 ) \, dx \, dy = 2 \iint_{\R^2_+} ( \QQ_v^e \times \pt_x \QQ_v^e) \cdot \pt_y \QQ_v^e \, dx \, dy.
\end{equation}
Since $|\QQ_v^2| \leq 1$ on $\R^2_+$ by the maximum principle, we can apply the Cauchy-Schwarz inequality in $\R^3$ to find
$$
2|( \QQ_v^e \times \pt_x \QQ_v^e) \cdot \pt_y \QQ_v^e| \leq (|\pt_x \QQ_v^e|^2 + |\pt_y \QQ_v^e|^2).
$$
Thus we conclude that $|\pt_x \QQ_v^e| = |\pt_y \QQ_v^e| \equiv 0$ from \eqref{eq:poho} if $|v| \geq 1$ (where a little refinement is needed in the limiting case $|v|=1$.) Hence $\QQ_v^e$ and therefore $\QQ_v$ must be constant. Again, we refer the reader to \cite{FrSc-18} for details of this argument.

The completes our sketch of the proof of Theorem \ref{thm:Qv_class}.
\end{proof}

\subsection{Spectral Analysis of the Linearised Operator}

We shall now briefly review the results in \cite{LeSc-18} on the spectrum of the linearised operator around traveling solitary waves for (HWM) in the energy-critical dimension $d=1$. More specifically, by following \cite{LeSc-18}, we consider the static case with vanishing velocity $v=0$, i.\,e., half-harmonic maps from $\R$ to $\Sb^2$. Moreover, we focus on the special case when the profile $\QQ \in \dot{H}^{\frac 1 2}(\R; \Sb^2)$ is given by Blaschke factor of degree $m \geq 1$ having identical factors, i.\,e., the zeros $z_1 = \ldots = z_m \in \C_+$ in Theorem \ref{thm:Qv_class} coincide. (For $m=1$ this is no loss of generality, whereas for $m > 1$ this assumption is non-trivial.) Henceforth we assume that $z_k = +\im$ for all $k=1, \ldots, m$. Thus we deal with profiles of the form
\begin{equation} \label{eq:pure_Blaschke}
\QQ(x) =  \QQ_m(x) = \left ( \mathrm{Re}\,  Q_m(x), \mathrm{Im} \, Q_m(x), 0 \right ) \quad \mbox{with} \quad Q_m(z) = \left ( \frac{z-\im}{z+\im} \right )^m.
\end{equation}
To express perturbations around $\QQ : \R \to \Sb^2$, we use the frame $\{ \mathbf{e}, J \mathbf{e} \}$ for vectors in the tangent space $T_\QQ \Sb^2$, where $\mathbf{e}=(0,0,1)$ and $J \mathbf{e} = \QQ \times \mathbf{e}$. In this frame, the linearized equation (HWM) around $\QQ$ is  found to be
\begin{equation} \label{eq:linearized_hwm}
\pt_t \left [ \begin{array}{c} h_1 \\ h_2 \end{array} \right ]  = J L   \left [ \begin{array}{c} h_1 \\ h_2 \end{array} \right ] + O(\mathbf{h}^2) 
\end{equation}
for the perturbation $\mathbf{h} = h_1 \mathbf{e} + h_2 J \mathbf{e}$. Here $O(\mathbf{h}^2)$ stands for quadratic terms in $\mathbf{h}$ and $|\nabla| \mathbf{h}$. The linearised operator is found to be
\begin{equation}
J L = \left [ \begin{array}{cc} 0 & -1 \\ 1 & 0 \end{array} \right ]  \left [ \begin{array}{cc} L_+ & 0 \\ 0 & L_- \end{array} \right ] = \left [ \begin{array}{cc} 0 & - L_- \\ L_+ & 0  \end{array} \right ].
\end{equation}
Here $L_+$ and $L_-$ are the scalar operators given by
\begin{equation}
L_+ = |\nabla| - | \pt_x \QQ_m | = |\nabla| - \frac{2m}{1+x^2} \quad \mbox{and} \quad L_- = L_+  + R ,
\end{equation}
and $R$ denotes the integral operator of the form
\begin{equation}
(R f)(x) = \frac{1}{2 \pi} \int_{\R} \frac{|\QQ_m(x)- \QQ_m(y)|^2}{|x-y|^2} f(y) \, dy. 
\end{equation}

Let us first consider the {\em nullspaces} of the operators $L_+$ and $L_-$, which are defined as
$$
\mathcal{N}(L_+) = \{ f \in \dot{H}^{\frac 1 2}(\R) \mid  L_+ f = 0 \} \quad \mbox{and} \quad \mathcal{N}(L_-) = \{ f \in \dot{H}^{\frac 1 2}(\R) \mid L_- f =0 \}.
$$
Note that also (non-trivial) constant functions belong to $\dot{H}^{\frac 1 2}(\R)$. Elements $f \in \mathcal{N}(L_\pm) \setminus L^2(\R)$ will be called {\em resonances} of $L_+$ and $L_-$, respectively.  From \cite{LeSc-18} we recall the following complete description of the nullspaces.

\begin{theorem} \label{thm:spec1}
For any degree $m \geq 1$, it holds that
$$
\mathcal{N}(L_+) = \mathrm{span} \, \{ f_m , g_m \} \quad \mbox{and} \quad \mathcal{N}(L_-) = \mathrm{span} \, \{ 1, f_1, \ldots, f_m, g_1, \ldots, g_m \},
$$ 
where the functions $\{ f_k \}_{k=1}^m$ and $\{g_k \}_{k=1}^m$ are known in closed form. In particular, we have that 
$$
\dim \mathcal{N}(L_+) = 2 \quad \mbox{and} \quad \dim \mathcal{N}(L_-) = 2m+1
$$
and hence {\bf nondegeneracy} holds for the linearised operator around half-harmonic maps $\QQ= \QQ_m \in \dot{H}^{\frac 1 2}(\R;\Sb^2)$ for any degree $m \geq 1$.

Moreover, the operators $L_+$ and $L_-$ have a common resonance $\varphi \in \dot{H}^{\frac 1 2}(\R) \setminus L^2(\R)$, which is given by 
$$
\varphi = \begin{dcases*} f_m & if $m$ is odd, \\ g_m & if $m$ is even. \end{dcases*} 
$$
\end{theorem}

\begin{remarks*}
{\em 1. The idea of the proof will be sketched below.

2. The notion of nondegeneracy means that all the elements in the nullspaces of $L_+$ and $L_-$ are entirely generated by the continuous symmetries of the family of half-harmonic maps, i.\,e., due to rotations, translations and changes of the $m$ zeros in the corresponding Blaschke product; see \cite{LeSc-18} for details.

3. For the special case $m=1$, the nondegeneracy of half-harmonic maps was also shown in \cite{SiWeZh-18}.   } 
\end{remarks*}

The spectral analysis of the operators $L_+$ and $L_-$ can be extended much further, as done in \cite{LeSc-18}. To this end, we recall that the point spectrum of the unbounded self-adjoint operators $L_+$ and $L_-$ acting on $L^2(\R)$ is defined as the set of eigenvalues in $L^2(\R)$. In \cite{LeSc-18}, the following collection of results is established.

\begin{theorem} \label{thm:spec2}
For any degree $m \geq 1$, we have the following.
\begin{itemize}
\item[(i)] {\bf $L^2$-Eigenvalues of $L_+$:} The point spectrum of $L_+$  consists of exactly $2m$ eigenvalues, i.\,e., 
$$
\sigma_\mathrm{p}(L_+) = \left \{ E_{0}, E_1, \ldots, E_{2m-1} \right  \}.
$$
Moreover, each eigenvalue $E_k$ is {\bf simple} and we have the inequalities
$$
E_0 < E_1 < \ldots < E_{2m-2} < E_{2m-1} = 0.
$$

\item[(ii)] {\bf $L^2$-Eigenvalues of $L_-$:} The point spectrum of $L_-$ only contains zero, i.\,e.,
$$
\sigma_{\mathrm{p}}(L_-) = \{ 0 \}.
$$
Moreover, the eigenvalue $E=0$ is exactly $2m$-fold degenerate.

\item[(iii)] Both operators $L_+$ and $L_-$ have {\bf no embedded} eigenvalue $E > 0$ inside  $\sigma_{\mathrm{ess}}(L_+) = \sigma_{\mathrm{ess}}(L_-) = [0, \infty)$.
\end{itemize}
\end{theorem}

\begin{remarks*}
{\em 1. A particular corollary of the preceding result is a {\em coercivity estimate} for $L_+$ and $L_-$ under certain (natural) orthogonality conditions. Such a result is typically needed in the modulational analysis of solutions close to $\QQ$.

2. It can be shown that the {\em purely absolutely continuous spectrum} $\sigma_{\mathrm{ac}}(L_+) = \sigma_{\mathrm{ac}}(L_-)$ equals $[0, \infty)$. For details on an explicit unitary transform establishing this fact,  see \cite{LeSc-18}.

3. For the spectral analysis of the matrix operator $\mathcal{L} = JL$, we refer again to \cite{LeSc-18}.}
\end{remarks*}

The proofs of Theorems \ref{thm:spec1} and \ref{thm:spec2} both essentially depend on the use of the stereographic projection $\Pi$ from the unit circle $\Sb$ to the projective real line $\hat{\R} = \R \cup \{ \infty \}$. To briefly sketch the main ideas, we note that,  by applying $\Pi$, we can recast the spectral analysis for the operators
$$
L_+= |\nabla| - \frac{2m}{1+x^2} \quad \mbox{and} \quad L_-= |\nabla| - \frac{2m}{1+x^2} + R
$$
in terms of the unbounded operators given by
$$
J = (1-\sin \theta) ( |\nabla|_{\Sb} - m \mathds{1}) \quad \mbox{and} \quad H = (1- \sin \theta) ( |\nabla|_{\Sb} - m \mathds{1} + \widetilde{R} )
$$
acting on $L^2(\Sb)$. Here $|\nabla|_{\Sb}$ is the square root of the Laplacian on $\Sb$ and $\widetilde{R}$ is some integral operator. Clearly, the operators $J \neq J^*$ and $H\neq H^*$ fail to be self-adjoint on $L^2(\Sb)$. Hence their detailed spectral analysis seems to be a hopeless enterprise at first sight. However, the benefit of this approach is that $J$ and $H$ are both seen to be {\bf Jacobi operators}. By this, we mean that the corresponding (infinite) matrices have a tridiagonal structure with respect to the standard Fourier basis $\{ e^{\im k \theta} \}_{k \in \mathbb{Z}}$ of $L^2(\Sb)$. For instance, the matrix for $J$ is of the form
$$
[J_{kl}] = \left [ \begin{array}{cccccc} \ddots & \ddots & \ddots &  & &  0 \\ & a_{n} & b_{n} & c_{n} \\ & & a_{n+1} & b_{n+1} & c_{n+1}   \\ 0 & &  & \ddots & \ddots & \ddots  \end{array} \right ] 
$$
with certain sequences $(a_n)$, $(b_n)$, $(c_n)$  and $n \in \mathbb{Z}$. Roughly speaking, the degree $m \in \mathbb{N}$ of the half-harmonic maps $\QQ = \QQ_m$ plays a central role by splitting the frequencies $k \in \mathbb{Z}$ on the unit circle with $|k| \leq m$ and $|k| \geq m$. More precisely, the analysis of the Jacobi operators given by $J$ and $H$ then shows the following. 
\begin{itemize}
\item The {\em bound states} of $L_+$ and $L_-$ are determined via the actions of $J$ and $H$ on the  $2m+1$-dimensional subspace $\mathrm{span} \, \{ e^{\im k \theta} : |k| \leq m \}$ in $L^2(\Sb)$. 
\item The {\em scattering states} of $L_+$ and $L_-$  can be analysed in detail via the action of $J$ and $H$ on the infinite-dimensional subspace  $\mathrm{span} \, \{ e^{\im k \theta} : |k| \geq m \}$ in $L^2(\Sb)$.
\end{itemize}
A distinguished part will be played by the two-dimensional subspace spanned by $\{ e^{\im k \theta} : k = \pm m \}$, which lies at the interface between the bound states and scattering states (i.\,e.~the continuous spectrum), and which yields  the two linearly independent solutions of $L_+ \varphi = L_-\varphi = 0$ given by the $L^2$-zero mode and the zero-energy resonance.

\section{Lax Pair Structure and Rational Solutions} \label{sec:Lax}

In the recent work \cite{GeLe-18}, a Lax pair was found for one-dimensional (HWM). This indicates that some sort of complete integrability is present in this case. To formulate this result (and some of its applications), we need to introduce some notation as follows. Let $\{ \sigma_1, \sigma_2, \sigma_3 \}$ denote the standard Pauli matrices. For notational convenience, we make use of the Pauli vector $\bm{\sigma}=(\sigma_1, \sigma_2, \sigma_3)$ in the following. Given a vector $\mathbf{x} \in \R^3$, we define the complex $2 \times 2$-matrix by setting
$$
\Xsp =  \mathbf{x} \cdot \bm{\sigma} = \sum_{j=1}^3 x_j \sigma_j = \left [ \begin{array}{cc} x_3 & x_1 - \im x_2 \\ x_1 + \im x_2 & -x_3 \end{array} \right ] .
$$
Clearly, we have that $\Xsp = \Xsp^*$ is Hermitian matrix with zero trace $\mathrm{Tr} \,\Xsp = 0$. From elementary algebra for the Pauli matrices we get the matrix identity
\begin{equation} \label{eq:nice}
( \mathbf{x} \cdot \bm{\sigma})(\mathbf{y} \cdot \bm{\sigma}) = (\mathbf{x} \cdot \mathbf{y}) \mathds{1} + \im ( \mathbf{x} \times \mathbf{y}) \cdot \bm{\sigma}.
\end{equation}
for all vectors $\mathbf{x}, \mathbf{y} \in \R^3$, where $\mathds{1}$ denotes the $2 \times 2$-unit matrix. As a consequence of this identity, we find that (HWM) can be written in the form
\begin{equation} 
\pt_t \Usp = -\frac{\im}{2} \left [ \Usp, \Dh \Usp  \right ],
\end{equation}
where $[A,B] = A B-BA$ is the commutator of matrices in $\C^{2 \times 2}$ and $\Usp = \uu \cdot \bm{\sigma}$. In abstract terms, the $\Sb^2$-valued map $\uu$ can be expressed in terms of the map $\Usp$ that takes values in the real Lie algebra $\mathfrak{su}(2)$ spanned by the Pauli matrices. Note that $\Usp^2 = |\uu|^2 \mathds{1} = \mathds{1}$ holds  as an immediate consequence of \eqref{eq:nice}.

Given a map $\uu : [0,T) \times \R \to \Sb^2$ (not necessarily a solution of (HWM)), we set $\Usp= \uu \cdot \bm{\sigma}$ and define the operators
\begin{equation} \label{def:LB}
L _\uu= [H, \Usp] \quad \mbox{and} \quad B_\uu = -\frac{\im}{2} \left (  \Usp \circ \Dh + \Dh \circ \Usp \right ) + \frac{\im}{2} \Dh \Usp
\end{equation}
acting on $L^2(\R; \C^2)$. Here $H$ denotes {\em Hilbert transform} on  $L^2(\R)$ defined in Fourier space by $\widehat{(H f)}(\xi) =-\im (\mathrm{sgn} \, \xi) \widehat{f}(\xi)$. For $\C^2$-valued functions $f \in L^2(\R; \C^2)$, we let $H$ act componentwise via $Hf = (H f_1, Hf_2)$ with a slight abuse of notation.

In the definition of $L$ and $B$, the matrix-valued functions $\Usp$ and $|\nabla|\Usp$ are understood as multiplication operators on $L^2(\R; \C^2)$. By the skew-symmetry of the Hilbert transform $H=-H^*$, we deduce the formal properties
$$
L_\uu =(L_\uu)^* \quad \mbox{and} \quad B_\uu=-(B_\uu)^*.
$$
Indeed, as we will see below, the operator $L_\uu$ is  of Hilbert-Schmidt class (and hence compact) if and only if $\uu$ has finite energy $E[\uu] < \infty$. Of course, the operator $B_\uu$ is always an unbounded operator on $L^2(\R; \C^2)$ due to the presence of the pseudo-differential operator $|\nabla|$. 

From \cite{GeLe-18} we recall the following result.

\begin{theorem}[Existence of a Lax Pair] \label{thm:Lax}
Let $\uu : [0,T) \times \R \to \Sb^2$ be a sufficiently regular solution of (HWM) in $d=1$ dimension. Then the following Lax equation holds true:
$$
\frac{d}{dt} L_\uu  = [B_\uu , L_\uu],
$$
where the operators $L_\uu$ and $B_\uu$ are defined in \eqref{def:LB} above with $\Usp=\uu \cdot \bm{\sigma}=\sum_{j=1}^3 u_j \sigma_j$.
\end{theorem}

The proof of Theorem \ref{thm:Lax} follows from a tricky calculation involving the fact that $|\nabla| = H \frac{d}{dx}$ holds combined with {\em Cotlar's product identity} for the Hilbert transform:
\begin{equation}
H(fg) = (Hf) g + f(Hg) + H((Hf)(Hg))
\end{equation}
for all $f,g \in \mathcal{S}(\R)$. Also, the normalisation condition $\Usp^2 = |\uu|^2 \mathds{1} = \mathds{1}$ enters in essential way into the arguments. We refer to \cite{GeLe-18} for the complete proof of Theorem \ref{thm:Lax}.

From the Lax identity above, we can obtain two fundamental corollaries about (HWM) in one space dimension:
\begin{enumerate}
\item Infinite number of conservation laws.
\item Rational initial data give rise to rational solutions.
\end{enumerate}

To elaborate on these claims, we note that the spectrum of $L_\uu$ stays constant in time thanks to the Lax equation in Theorem \ref{thm:Lax}. As an immediate consequence, we obtain the following result on an infinite number of (formal) conservation laws for the one-dimensional (HWM).

\begin{corollary} \label{cor:L_trace}
For $\uu : [0,T) \times \R \to \Sb^2$ as above, we have the (formal) conservation laws
$$
\mathrm{Tr} ( |L_\uu|^p ) = \mathrm{const}. 
$$
for any $0 < p < \infty$.
\end{corollary}

Let us make some few comments on this result. If we take $p=2$ above, we obtain the conservation of energy again. Indeed, the integral kernel of the commutator $L_\uu = [H, \Usp]$ is found to be
$$
K(x,y) = \frac{1}{\pi} \frac{\Usp(x)- \Usp(y)}{x-y} \in \C^{2 \times 2},
$$
where we omit the $t$-dependence for notational convenience. Using that $\mathrm{Tr}_{\C^2} (\Usp^2) = \mathrm{Tr}_{\C^2}\, (\mathds{1} |u|^2) = 2 |\uu|^2$, we find that the square of the Hilbert-Schmidt norm of $L_\uu$ is
$$
\mathrm{Tr} \, ( |L_\uu|^2) = \frac{2}{\pi^2} \iint_{\R\times \R} \frac{|\uu(x)-\uu(y)|^2}{|x-y|^2} \, dx \,dy = \frac{8}{\pi} E[\uu].
$$
In particular, we see that $L_\uu$ is of Hilbert-Schmidt class if and only if $\uu \in \dot{H}^{\frac 1 2}(\R; \Sb^2)$.

For general $0 < p < \infty$, we can extend the previous observation by exploiting the norm equivalence of {\em Schatten norms} of the operator $L_\uu$ with certain {\em homogeneous Besov norms} of function $\uu$. More precisely, let us recall that the family of Schatten norms is given by  
\begin{equation}
\| L_\uu \|_{\mathfrak{A}_p} = (\mathrm{Tr} \, ( |L_\uu|^p) )^{\frac{1}{p}} \quad \mbox{for} \quad 0 < p < \infty.
\end{equation}
For $p \geq 1$, this in fact a norm, whereas for $0 < p < 1$ we only obtain a quasi-norm. By adapting a classical result due to V.~Peller on Hankel operators (see e.\,g.~\cite{Pe-04} and  by using the splitting $L^2(\R) = L_+^2(\R) \oplus L_-^2(\R)$ into positive and negative frequencies, we obtain the equivalence
\begin{equation}
\| L_\uu \|_{\mathfrak{A}_p}  \sim_p \| \uu \|_{\dot{B}_{p,p}^{1/p}}  \quad \mbox{for} \quad  0 < p < \infty.
\end{equation}
Note that $\| \cdot  \|_{\dot{B}_{p,p}^{1/p}}$ is a semi-norm if and only if $p \geq 1$, where it is only a quasi semi-norm when $0 < p < 1$. By Corollary \ref{cor:L_trace} and classical facts about Besov space, we deduce the a-priori bound
\begin{equation}
\sup_{t \in [0,T)} \left ( \int_{\R} |\pt_x \uu(t,x)| \, dx \right ) \lesssim \sup_{t \in [0,T)} \| \pt_x \uu(t) \|_{B^0_{1,1}} \lesssim \| \uu(0) \|_{\dot{B}^1_{1,1}}.
\end{equation}
In geometric terms, this a-priori estimate implies that the {\em length of the curve} parametrised by $\uu (t) : \R \to \Sb^2$ is uniformly bounded for all $t \in [0,T)$. However, it remains a very interesting open problem to show that the conservation laws obtained from $L_\uu$ are strong enough to deduce global well-posedness result for the one-dimensional (HWM). Furthermore, it is possible that although a smooth solution $\uu=\uu(t,x)$ exists for all times $t \geq 0$, higher Sobolev norms $\| \uu(t) \|_{\dot{H}^{\frac{1}{2}+ \eps}}$ could growth drastically in $t$, showing some sort of turbulent behaviour. Such a turbulence phenomenon was recently constructed in \cite{GeLePoRa-18} for the focusing cubic half-wave equation in $\R$, exploiting the existence of small traveling solitary waves. It may be conjectured that an analogous behaviour can happen for (HWM) building upon the existence of small traveling solitary waves, see Section \ref{sec:solitary} above.

Another main corollary from Theorem \ref{thm:Lax} follows from an adaptation of {\em Kronecker's theorem} for Hankel operators, which says that $L_\uu$ has finite rank if and only if $\uu=\uu(t,x)$ is a rational function of $x$. Indeed, we have the following.

\begin{corollary}
Let $\uu : [0,T) \times \R \to \Sb^2$ solve (HWM) in $d=1$ dimension. If the initial datum $\uu(0,x)$ is a rational function of $x$, then $\uu(t,x)$ is a rational function of $x$ for all $t \in [0,T)$.
\end{corollary}

The idea of the proof involves showing that $L_\uu=[H, \Usp]$ is of finite rank if and only if $\uu(t,x)$ is rational in $x$. In fact, we can relate the rank of $L_\uu$ with the number of poles of $x \mapsto \uu(t,x)$ in a precise way:
$$
\mathrm{rank} \, L_\uu = \mbox{number of poles of $\uu(x)$ in $\C \setminus \R$}.
$$
By the Lax type evolution stated in Theorem \ref{thm:Lax}, the rank of $L_\uu$ stays constant in time and hence the result follows. Again, it is a very interesting open problem to show global-in-time existence in the restricted class of rational solutions to (HWM). Even in this setting, it not clear yet how to exploit the conservation laws to deduce global-in-time existence.

Note that a special subclass of rational solutions is given by the traveling solitary waves $\uu(t,x) = \mathbf{Q}_v(x-vt)$ discussed in Section \ref{sec:solitary} above. Furthermore, in the class of rational solutions for (HWM), there exist other interesting explicit solutions such as {\em time-periodic solutions}. For instance, the function (which is rational in $x$) given by
\begin{equation}
\uu_{\mathrm{per}}(t,x) = \left ( \cos \left ( \frac{t}{\sqrt{2}} \right) \frac{2x^2}{x^4+1}, \sin \left ( \frac{t}{\sqrt{2}} \right ) \frac{2x^2}{x^4+1}, \frac{x^4-1}{x^4+1} \right )
\end{equation}
is a time-periodic solution of the one-dimensional (HWM). At each time $t$ fixed, the map $x \mapsto \uu_{\mathrm{per}}(t,x)$ parametrises a half-circle on $\Sb^2$ starting from $\uu_{\mathrm{per}}(t,-\infty) = (0,0,1)$ and going back to $\uu_{\mathrm{per}}(t,+\infty)=(0,0,1)$. Its time-dependence corresponds to a rotation of this half-circle around the $x_3$-axis with constant angular velocity $\omega = 1/\sqrt{2}$. It is an interesting open problem to study the stability properties of $\uu_{\mathrm{per}}(t,x)$ under rational (or more general) perturbations.

\section{Periodic Case and Hyperbolic Target $\mathbb{H}^2$} \label{sec:h2}

A natural variant of (HWM) arises in the one-dimensional periodic setting when  $\R$ is replaced by $\mathbb{T} = \R / 2 \pi \mathbb{Z}$. That is, we consider maps $\uu : [0,T) \times \mathbb{T} \to \Sb^2$ solving
\begin{equation}
\pt_t \uu = \uu \times |\nabla| \uu
\end{equation}
where $|\nabla| f = \sum_{n \in \mathbb{Z}} |n| \widehat{f}_n e^{\im n \theta}$ now denotes the square root of the Laplacian on $\mathbb{T}$. By means of the stereographic projection $\Pi : \Sb \to \R \cup \{ \infty \}$ it can be shown that the traveling solitary wave profiles $\QQ_v : \R \to \Sb^2$ on the real line give rise to traveling solitary waves profile $\tilde{\QQ}_v : \mathbb{T} \to \Sb^2$ for (HWM) on $\mathbb{T}$ and vice versa, by the simple relation $\tilde{\QQ}_v = \QQ_v \circ \Pi$. This fact also reflects the conformal invariance of the energy functional $E[\uu]$ as sketched above; see \cite{LeSc-18}. In particular, the  classification result in Theorem \ref{thm:Qv_class} applies mutatis mutandis to the one-dimensional periodic (HWM). Furthermore, the existence and the formula for the Lax pair $(L_\uu, B_\uu)$ both carry over to periodic case in a direct fashion; see \cite{GeLe-18}.

A further natural variation of (HWM) occurs when the compact target $\Sb^2$ is replaced by the non-compact K\"ahler manifold $\mathbb{H}^2$, i.\,e., the two-dimensional hyperbolic plane. Again, we refer to \cite{GeLe-18} for more details and the construction of the Lax pair for the target $\mathbb{H}^2$.

\section{Conclusion and (Some) Open Problems}

In this short note, we have given a brief overview on the current state of affairs concerning the half-wave maps equation (HWM). This geometric evolution equation exhibits a set of analytically intriguing properties. In the energy-critical case of $d=1$ space dimensions, we obtain: complete classification of traveling solitary waves, detailed spectral theory of linearised operators, time-periodic solutions, the existence of a Lax pair, infinitely many conservation laws involving the theory of Hankel operators, invariance of rational solutions under the flow, etc. Yet, it is fair to say that we are still at the beginning of understanding the dynamical properties of solutions. In the author's opinion, the following three circles of open problems are of main interest for a future study of the energy-critical (HWM):
\begin{enumerate}
\item Global well-posedness for (HWM) in energy space (or counterexamples by singularity formation).
\item Stability (or instability) of traveling solitary waves or other special solutions (e.\,g.~the time-periodic rational solutions).
\item Develop a notion of complete integrability building on the Lax pair.
\end{enumerate}
Any progress concerning (1)--(3) would be highly desirable.

\subsubsection*{Acknowledgments}
The author is grateful to Patrick G\'erard and Armin Schikorra for fruitful collaborations on (HWM).


\begin{thebibliography}{99}

\bibitem{DaLioRi-11} F.~Da~Lio and T.~Rivi{\`e}re, Sub-criticality of non-local Schr\"odinger systems with antisymmetric potentials and applications to half-harmonic maps, {\em Advances in Mathematics} 227 (2011), 1300 -- 1348.

\bibitem{GeLe-18} P. G{\'e}rard and E. Lenzmann, A Lax pair structure for the half-wave maps equation, {\em Lett. Math. Phys.} 108 (2018), no. 7, 1635--1648.

\bibitem{GeLePoRa-18} P. G{\'e}rard, E. Lenzmann, O. Pocovnicu, and P. Rapha{\"e}l, A two-soliton with transient turbulent regime for the cubic half-wave equation on the real line, {\em Annals of PDE} 4 (2018), no.1, 166pp.

\bibitem{FrSc-15} A. Fraser and R. Schoen, Uniqueness theorems for free boundary minimal disks in space forms, {\em Int. Math. Res. Not. IMRN} 17 (2015), 8268--8274.

\bibitem{KrSi-18} J. Krieger and Y. Sire,  Small data global regularity for half-wave maps, {\em Anal. PDE} 11 (2018), no. 3, 661--682.

\bibitem{LeSc-18} E. Lenzmann and A. Schikorra, On energy-critical half-wave maps into $\mathbb{S}^2$, {\em Invent. Math.} 213 (2018), no. 1, 1--82.

\bibitem{MiSi-15} V. Millot and Y. Sire, On a fractional Ginzburg-Landau equation and 1/2-harmonic maps into sphere, {\em Arch. Rat. Mech. Anal.} 215 (2015), no. 1, 125--210.

\bibitem{MiPi-04} P. Mironescu and A. Pisante, A variational problem with lack of compactness for {$H^{1/2}(S^1;S^1)$} maps of prescribed degree, {\em J. Funct. Anal.} 217 (2004), no. 2, 249--279.

\bibitem{Pe-04} V. Peller, Hankel operators and their applications, Springer, New York (2003).

\bibitem{PuGu-13} X. Pu and B. Guo, Well-posedness for the fractional Landau-Lifshitz equation without Gilbert damping, {\em Calc. Var. PDE} 46 (2013), no. 3--4, 441--460.  

\bibitem{SiWeZh-18} Y. Sire, J. Wei, and Y. Zheng, Nondegeneracy of half-harmonic maps from $\mathbb{R}$ into $\mathbb{S}^1$, {\em Proc. Amer. Math. Soc.} 146 (2018), no. 12, 5263--5268. 
 
\bibitem{SiWeZh-18b} Y. Sire, J. Wei, and Y. Zhang, Infinite time blow-up for half-harmonic map flow from $\mathbb{R}$ into $\mathbb{S}^1$, {\em arXiv:1711.05387}.


\bibitem{ZhSt-15} T.~Zhou and M.~Stone, Solitons in a continuous classical {Haldane--Shastry} spin chain {\em Phys. Lett. A} 379 (2015), 2817--2825.


\end{thebibliography}
\end{document}